\def\draft{n}
\documentclass[12pt]{amsart}
\usepackage[headings]{fullpage}
\usepackage{amssymb,epic,eepic,epsfig}%,amsbsy,amsmath,bbold}
\usepackage{epstopdf}
\usepackage{graphicx}
\usepackage{texdraw}
\usepackage{url}
\usepackage[bookmarks=true,%
    colorlinks=true,%
    linkcolor=blue,%
    citecolor=blue,%
    filecolor=blue,%
    menucolor=blue,%
    urlcolor=blue,%
    breaklinks=true]{hyperref}

\theoremstyle{definition}

\def\printname#1{
        \if\draft y
                \smash{\makebox[0pt]{\hspace{-0.5in}
                        \raisebox{8pt}{\tt\tiny #1}}}
        \fi
}

\newlength{\standardunitlength}
\catcode`\@=11
\long\def\@makecaption#1#2{%
     \vskip 10pt

\setbox\@tempboxa\hbox{%\ifvoid\tinybox\else\box\tinybox\fi
       \small\sf{\bfcaptionfont #1. }\ignorespaces #2}%
     \ifdim \wd\@tempboxa >\captionwidth {%
         \rightskip=\@captionmargin\leftskip=\@captionmargin
         \unhbox\@tempboxa\par}%
       \else
         \hbox to\hsize{\hfil\box\@tempboxa\hfil}%
     \fi}
\font\bfcaptionfont=cmssbx10 scaled \magstephalf
\newdimen\@captionmargin\@captionmargin=2\parindent
\newdimen\captionwidth\captionwidth=\hsize
\catcode`\@=12

\def\lbl#1{\label{#1}\printname{#1}}

\def\BQ{\mathbb Q}

\def\BC{\mathbb C}

\def\calT{\mathcal T}

\def\ga{\gamma}

\def\la{\langle}
\def\ra{\rangle}

\def\calB{\mathcal{B}}

\def\Vol{\mathrm{Vol}}

\def\Li{\mathrm{Li}}
\def\calB{\mathcal{B}}

\newcommand{\mb}{\mathbf}

\begin{document}

%%%%%%%%%%%%%%%%%%%%%%{page1}

\title[Exact computation of the $n$-loop invariants of knots]{
Exact computation of the $n$-loop invariants of knots}

\author{Stavros Garoufalidis}
\address{School of Mathematics \\
         Georgia Institute of Technology \\
         Atlanta, GA 30332-0160, USA \newline
         {\tt \url{http://www.math.gatech.edu/~stavros }}}
\email{stavros@math.gatech.edu}
\author{Eric Sabo}
\address{School of Mathematics \\
         Georgia Institute of Technology \\
         Atlanta, GA 30332-0160, USA \newline
         {\tt \url{http://www.math.gatech.edu/users/esabo3 }}}
\email{esabo3@gatech.edu}
\author{Shane Scott}
\address{School of Mathematics \\
         Georgia Institute of Technology \\
         Atlanta, GA 30332-0160, USA \newline
         {\tt \url{http://www.math.gatech.edu/users/sscott42 }}}
\email{scottsha@gatech.edu}
\thanks{
S.G. was supported in part by NSF.
1991 {\em Mathematics Classification.} Primary 57N10. Secondary 57M25.
\newline
{\em Key words and phrases: knots, Jones polynomial, Kashaev invariant,
Volume Conjecture, volume, hyperbolic geometry, ideal triangulations,
shapes, Neumann-Zagier data, 1-loop, n-loop, formal Gaussian integration,
Feynman diagrams.}
}

\date{March 24, 2015}%\today}

%\dedicatory{\large{\bf Private notes. Please do not
%distribute under any circumstances!}}

\begin{abstract}
The loop invariants of Dimofte-Garoufalidis is a formal power series with
arithmetically interesting coefficients that conjecturally appears
in the asymptotics of the Kashaev invariant of a knot to all orders in $1/N$.
We develop methods implemented in \texttt{SnapPy} that compute the first 6 
coefficients of the formal power series of a knot. We give examples that
illustrate our method and its results.
\end{abstract}

\maketitle

\tableofcontents

%%%%%%%%%%%%%%%%%%%%%%%%%%%%%%%%%%%%%%%%%%%%%%%%%%%%%%%%%%%%%%%%%%%%%%%%%%%%
%%%%%%%%%%%%%%%%%%%%%%%%%%%%%%%%%%%%%%%%%%%%%%%%%%%%%%%%%%%%%%%%%%%%%%%%%%%%

%%%% programs/sage/feynman.diagrams

\section{Introduction}
\lbl{sec.intro}

\subsection{The Volume Conjecture to all orders in $1/N$}
\lbl{sub.nloop}

The best known quantum invariant of a knot in 3-space is the Jones 
polynomial~\cite{Jones}. The Kashaev invariant $\la K \ra_N$ of a knot
$K$ (for $N=1,2,\dots$)~\cite{Kashaev95} coincides with the evaluation of the 
Jones polynomial of a knot and its parallels at complex roots of 
unity~\cite{MM}. The Volume Conjecture of Kashaev~\cite{K97} states 
that for a hyperbolic knot $K$, 
$$
\lim_{N \to \infty} \frac{1}{N} \log | \la K \ra_N| = \frac{\Vol(K)}{2\pi},
$$
where $\Vol(K)$ is the hyperbolic volume of $K$. An extension of the 
Volume Conjecture to all orders in $1/N$ was proposed independently by 
Gukov and the first author~\cite{Gu,Ga:arithmetic}. Namely, for every 
hyperbolic knot $K$ there exists a formal power series 
$\phi_K(\hbar) \in \BC[\![\hbar]\!]$ such that
\begin{equation}
\lbl{eq.expN}
\langle K \rangle_N \sim N^{3/2} e^{C_K N} \phi_K(2 \pi i/N),
\end{equation}
where $C_K$ is the complexified volume of $K$ divided by $2 \pi i$,
\begin{subequations}
\begin{align}
\lbl{eq.FK1}
\phi_K(\hbar) &=\tau_K^{-\frac{1}{2}} \,\, \phi^+_{K,1}(\hbar) ,\\
\lbl{eq.FK2}
\phi^+_K(\hbar) & \in 1 + \hbar F_K[\![\hbar]\!] ,\\
\lbl{eq.FK3}
\tau_K & \in F_K,
\end{align}
\end{subequations}
and $F_K$ is the trace field of $K$. 

\subsection{Ideal triangulations, shapes, and the loop invariants}
\lbl{sub.loop}

The left hand side of Equation~\eqref{eq.expN} is concretely
defined given a planar projection or an ideal triangulation of a knot,
and is typically given by a finite state-sum where the summand is a ratio
of quantum factorials. Examples of state-sum formulas for the Kashaev invariant
of the $4_1$, $5_2$ and $6_1$ knots are given in \cite[(2.2)-(2.4)]{K97}.

On the other hand, the power series $\phi_K(\hbar)$ that conjecturally appears 
in the right hand side is not an explicit function of the knot. 
Numerical computations of the Kashaev invariant were performed by Zagier and 
G., and using numerical interpolation and a variety of guessing methods, it 
was possible to recognize the first few coefficients of the power series 
$\phi_K$ for several knots~\cite{GZ1}.

\begin{figure}[htb]
\includegraphics[width=1.5in]{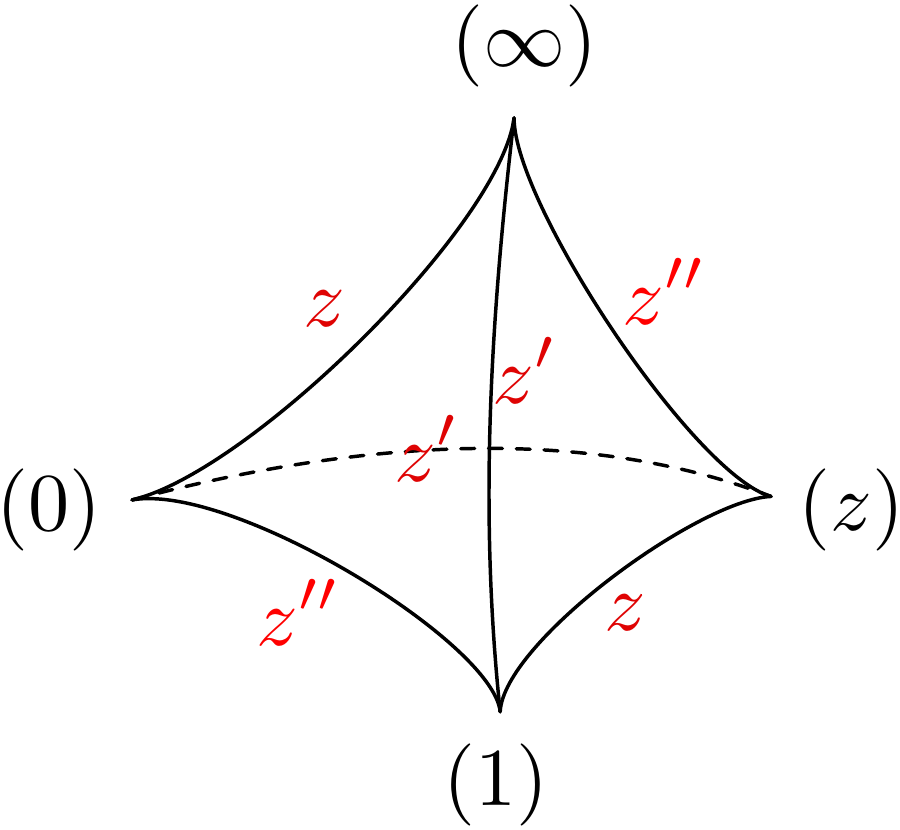}
\caption{An ideal tetrahedron and its shape assignment.}
\label{fig.tet}
\end{figure}

The main result of Dimofte-G.~\cite{DG} was the construction of a power series 
$\phi_{\gamma}(\hbar)$ that depends on an ideal triangulation of a knot 
complement. For a detailed discussion on ideal triangulations and their 
gluing equations, see~\cite{Th,NZ} and also \cite[Sec.1.2]{DG}.
Explicitly, an ideal triangulation $\calT$ with $N$ tetrahedra
gives rise to a vector $z=(z_1,\dots,z_N)$ of shapes that satisfy the 
Neumann-Zagier equations 
$$
\lbl{eq.NZ}
\prod_{j=1}^N z_j^{\mb A_{ij}} (z_j'')^{\mb B_{ij}}=(-1)^{\nu_i}
$$
for $i=1,\dots,N$ where $z_j''=1-1/z_j$. These equations are obtained as 
follows. Fix an edge of the ideal triangulation $\calT$, and set the
product of the shape parameters of all tetrahedra that go around the edge $e$
equal to $1$. If $z$ is the shape of a tetrahedron that contains the fixed 
edge, then its contribution to the above product is $z$ or $z'$ or $z''$ 
according to the convention of Figure \ref{fig.tet}. Finally, replace
$z'=-(z z'')^{-1}$ using the relation $z z' z'' = -1$. This gives rise
to the equations \eqref{eq.NZ}, one for each edge of $\calT$. Likewise,
there is an equation os the same type for each peripheral curve.

In the above equations, we have
removed one edge equation and replaced it with the meridian cusp equation.
Neumann-Zagier~\cite{NZ} prove that $(\mb A|\mb B)$ is the upper part of a 
symplectic matrix. It follows that $(\mb A|\mb B)$ has rank $N$ and 
$\mb A \mb B^T$ is symmetric, where $\mb B^T$ is the transpose of $\mb B$.
We will assume that $\mb B$ has nonvanishing determinant. Furthermore,
we will assume that our triangulation $\calT$ is such that there exists a 
solution to the gluing equations in $(\BC\setminus\{0,1\})^N$ that recovers
the complete hyperbolic structure of the hyperbolic knot $K$. In that
case, $z$ is a vector of algebraic numbers and $\BQ(z_1,\dots,z_N)$ is a 
number field (the shape field) which coincides with the invariant trace field
and with the trace field of the knot~\cite[Thm.2.2,2.4]{NR}. Finally, one
can choose a flatenning $(f,f'')$, that is an integer solution of the 
linear equation $\mb A f + \mb B f'' = \nu$. This determines a 
Neumann-Zagier datum $\ga=(A,B,\nu,f,f'',z)$, which in turn defines  
the power series $\phi_{\gamma}(\hbar)$. Of course, different ideal 
triangulations give rise to different Neumann-Zagier data, hence to 
potentially different formal power series $\phi_{\gamma}(\hbar)$. On the other 
hand, the left hand-side of Equation \eqref{eq.expN} depends only on the
hyperbolic knot $K$. Although the topological invariance of 
$\phi_{\gamma}(\hbar)$ is not known, from the computational point of view,
this gives an excellent consistency check of correctness of the code.

Equations~\eqref{eq.FK1}-\eqref{eq.FK3} are manifest by the definition of
$\phi_\gamma(\hbar)$. In~\cite{DG} it was shown that $\tau_\gamma$ is a 
topological invariant, defined up to a sign. We may call $\tau_\gamma$
the 1-loop invariant. If we write
$$
\phi^+_\gamma(\hbar)=\exp\left(\sum_{n=2}^\infty S_{\gamma,n}\hbar^{n-1}
\right),
$$
then $S_{\gamma,n}$ are the $n$-loop invariants of the $\ga$. In~\cite{DG}
it was conjectured that $S_{\gamma,2}$ is well-defined up to addition of
an integer multiple of $1/24$, and that $S_{\gamma,n}$ are topological
invariants for $n \geq 3$.
 
The definition of $\phi^+_\gamma(\hbar)$ is given explicitly by formal
Gaussian integration. It follows that $S_{\gamma,n}$ is a weighted sum
of a finite set of Feynman diagrams with Feynman loop number at most $n$. 
The Feynman rules were explained in
detail in~\cite[Sec.1.6-1.8]{DG} and the contributing Feynman diagrams for
$n=2$ and $n=3$ were explicitly drawn. For $n>3$, the number of Feynman 
diagrams gets large and drawings-by-hand is not advisable. 

For the benefit of the reader, we recall the Feynman rules 
from~\cite[Sec.1.6-1.8]{DG}. By connected \emph{Feynman diagram} $G$ we mean a 
connected multigraph, possibly with loops and multiple edges.
If $G$ is a Feynman diagram, its \emph{Feynman loop number} $L(G)$ is given by
$$
L(G)=|V_1(G)|+|V_2(G)|+b_1(G) 
$$
where $|V_k(G)|$ is the number of $k$-valent vertices of $G$ and 
$b_1(G)$ is the first betti number (also known as the number of holes) of $G$.
It is easy to see that a connected Feynman diagram with loop number
at most $n$ has at most $2n-2$ vertices and at most $n$ holes.
%%% see: feynGraphGeneration.tex and text2.tex
%%% These bounds are sharp. Indeed, if $n>2$ then the multigraph with a single 
%%% vertex and $n$ loops achieves the bound on the number of holes, and
%%% the trivalent graph with $2n-2$ vertices arranged in a cycle with every 
%%% other edge a 2-multiedge has achieves a loop number $n$.
Hence, there are finitely many Feynman diagrams of loop number at most $n$.

Fix a Neumann-Zagier datum $\ga=(\mb A, \mb B, \nu, f, f'', z)$ which
we assume is non-degenerate, that is the propagator (defined below) makes
sense. 
In each Feynman diagram $G$, the edges represent an $N\times N$ propagator
$$
\Pi=\hbar \left( -\mb B^{-1}\mb A+\mathrm{diag}(1/(1-z)) \right)^{-1}
$$
while each $k$-vertex comes with an $N$-vector of factors $\Gamma^{(k)}_i$,
$$
\Gamma^{(k)}_i = (-1)^k\sum_{p\,
=\,\alpha_k}^{\alpha_k +n-L(D)}\frac{\hbar^{p-1}(-1)^p B_p}{p!}\Li_{2-p-k}(z_i^{-1}) 
+ \begin{cases} -\tfrac12(\mb B^{-1}\nu)_i & k= 1 \\ \;\;0 & k \geq 2 
\end{cases}\,,
$$
where $\alpha_k = 1$ (resp., $0$) if $k=1,2$ (resp., $k\geq 3$).
Here $B_k$ is the $k$-th Bernoulli number ($B_1=-1/2$, $B_2=1/6$) and
$\Li_s(z)=\sum_{m=1}^\infty z^m/m^s \in \BQ(z)$ is the $s$-polylogarithm
function for $s$ a nonpositive integer.
The diagram $G$ is then evaluated by contracting the \emph{vertex factors} 
$\Gamma^{(k)}_i$ with 
propagators, multiplying by a standard {\em symmetry factor}, and taking the 
$\hbar^{n-1}$ part of the answer. In the end, $S_{\ga,n}$ is the sum of 
evaluated diagrams, plus an additional {\em vacuum} contribution
$$
\Gamma^{(0)} = \frac{B_n}{n!}\sum_{i=1}^N\Li_{2-n}(z_i^{-1}) 
+ \begin{cases} \tfrac18f\cdot\mb B^{-1}\mb Af & n=2 \\ 0 & n\geq 3 
\end{cases}\,.
$$

\subsection{Our code}
\lbl{sub.code}

Our goal is to give an exact computation for the $n$-loop invariants
for $n=1,\dots,6$ of a Neumann-Zagier datum of a \texttt{SnapPy}
triangulation. Our method is implemented in \texttt{SnapPy}.
We accomplished this in three steps.

\begin{itemize}
\item[(a)] 
We wrote a Python method 
\texttt{generate}\_\texttt{feynman}\_\texttt{diagrams.py}
that generates all Feynman diagrams that contribute to the $n$-loop 
invariant. The Feynman diagrams were generated by first generating trees,
and then adding to them multiple edges or loops. The number of such diagrams 
is shown in Table~\ref{t.diagrams}. 
Observe that if $G$ is a multigraph with corresponding simple graph $S(G)$ 
then $S(G)$ has at most $2n-2$ vertices and at most $n$ holes, and
$L(G)$ can be obtained from $S(G)$ by adding at most 
$n-L(S(G))+|V_1(S(G))|+|V_2(S(G))|$ edges. Thus, 
all Feynman diagrams with Feynman loop number at most $n$ can thus be 
generated by first generating all trees with at most $2n-2$ vertices then 
iteratively adding edges between pairs of vertices. Every edge added also 
adds an additional hole. If multigraph $G$ has more than $n-|V_1(G)|-|V_2(G)|$ 
holes it cannot be the subgraph of a Feynman diagram with Feynman loop 
number at most $n$. 
\begin{table}
\begin{center}
\begin{tabular}{|l||l|l|l|l|l|}\hline
 $n$ & $2$ & $3$ & $4$ & $5$ & $6$ \\ \hline
%$g_n$ & $6$ & $34$ & $291$ & $3369$ & $50058$ \\ \hline
$g_n$ & $6$ & $40$ & $331$ & $3700$ & $53758$ \\ \hline
\end{tabular}\vspace{.2cm}
\caption{The number $g_n$ of graphs that contribute to the $n$-loop 
invariant for $n=2,\dots,6$.} 
\lbl{t.diagrams}
\end{center}
\end{table}
\item[(b)]
We wrote a Python class \texttt{NeumannZagierDatum} which gives the 
Neumann-Zagier matrices and the exact value of the shape parameters
that recover the geometric representation of an ideal triangulation.
The exact computation of the shape parameters was done using the Ptolemy
module~\cite{Ga:ptolemy,snappy} and the numerical computation is already
implemented in \texttt{SnapPy}.
\item[(c)]
We wrote Python classes \texttt{nloop}\_\texttt{exact.py} 
and \texttt{nloop}\_\texttt{num.py} which given a Neumann-Zagier datum $\ga$
and a natural number $n=1,\dots,6$ computes $S_{\gamma,n}$ exactly (as 
an element of the trace field) or numerically to arbitrary precision.
\end{itemize}

To verify correctness of our code, we computed the $n$-loop invariants for
$n=1,\dots,5$ for different triangulations of each of a fixed knot, such 
as $5_2$, $(-2,3,7)$ pretzel, $6_1$, and $6_2$. In all cases the results 
agreed (up to a sign when $n=1$ and up to addition of $1/24$ times an integer 
when $n=2$). This illustrates both the topological invariance of the 
$n$-loop invariants, and the correctness of our code.

\subsection{Usage}
\lbl{sub.usage}

The essence of our code lies in two Python classes \texttt{NeumannZagierDatum}
and \texttt{nloop}. The former takes as input a manifold and generates the 
Neumann-Zagier datum $\ga=(A,B,\nu,f,f'',z)$, and the latter takes as input
Neumann-Zagier datum, an integer $n$, and a list of Feynman diagrams
and returns the $n$-th loop invariant $S_{\gamma,n}$.

The \texttt{NeumannZagierDatum} class has three optional arguments
\texttt{engine}, \texttt{verbose}, and \texttt{file\_name}, which are set to 
\texttt{None}, \texttt{False}, and \texttt{None}, respectively, by default. 
The \texttt{engine} variable is passed as an option into the Ptolemy module 
and controls the method in which solutions to the Ptolemy variety are found. 
The preferred value for this variable for manifolds in \texttt{CensusKnots} 
is \texttt{engine="magma"}, which refers to the \texttt{Sage} interface to 
the \texttt{Magma} Computational Algebra System \cite{magma}.
If Magma is not available, \texttt{engine="None"} will attempt to compute 
solutions of the Ptolemy variety using \texttt{Sage}. 
%This is known to be 
%currently unsuccessful with a handful of manifolds. (Note that in this case 
%the system simply hangs and no error message is produced.) 
Solutions for 
manifolds in \texttt{HTLinkExteriors} and \texttt{LinkExteriors} have been 
precomputed and are available with the Ptolemy module using 
\texttt{engine="retrieve"}, \cite{Ga:ptolemy}. This option requires an 
internet connection, but will automatically switch to recomputing locally 
if the download is unsuccessful. The output of the Ptolemy module including the
\texttt{retrieve} option are suppressed with \texttt{verbose=False} and are
displayed with \texttt{verbose=True}.

To utilize the \texttt{NeumannZagierDatum} class use a terminal to navigate to
the directory containing \texttt{nloop\_exact.py}, available at 
\cite{nloop-compute}, and load \texttt{Sage}. Once loaded, the class must 
first be initiated via\\
\indent \texttt{sage: attach('nloop\_exact.py')\\
\indent	     sage: M = Manifold('6\_1')\\
\indent	     sage: D = NeumannZagierDatum(M, engine="retrieve")}.\\
%Note that in this example we have used the optional argument
%\texttt{engine="retrieve"}, which will fail because the manifold $6_1$ is in
%\texttt{CensusKnots} and not \texttt{HTLinkExteriors} 
%or \texttt{LinkExteriors}.
%However, as noted above, such errors are automatically taken into account in
%\texttt{nloop\_exact.py}. 
To generate the Neumann-Zagier datum use\\
\indent \texttt{sage: D.generate\_nz\_data()}.\\
This will assign a Python list $[A,B,\nu,f,f'',z, \text{embedding}]$ 
consisting of
the Neumann-Zagier datum plus the embedding of the something in the
something to the class variable \texttt{nz}. If the optional argument
\texttt{file\_name} is used, this variable will be saved as a \texttt{Sage} 
object
file (*.sobj) in the current directory. To view the data simply use\\
\indent \texttt{sage: D.nz}.\\
The shape equations $z$ and field embedding may be computed separately
via\\
\indent \texttt{sage: D.exact\_shapes\_via\_ptolemy\_lifted()}\\
and\\
\indent \texttt{sage: D.compute\_ptolemy\_field\_and\_embedding()},\\
respectively.

Once the Neumann-Zagier datum has been computed, one may use it to 
compute the $n$-loop invariants $S_{\gamma,n}$. First, load the Feynman
diagrams you wish to use and choose an invariant you wish to calculate,\\
\indent \texttt{sage: n = 2}\\
\indent \texttt{sage: diagrams = load('6diagrams.sobj')}\\
\indent \texttt{sage: E = nloop(D.nz, n, diagrams)}.\\
Here, we have chosen to calculate $S_{\gamma,2}$ using Feynman
diagrams up to six loops. Note that the manifold \texttt{M} is not directly 
used when initiating the \texttt{nloop} class, as all the information about 
the manifold we
need is encoded in the Neumann-Zagier datum \texttt{D.nz}. To compute the
invariant use\\
\indent \texttt{sage: E.one\_loop()}\\
if $n = 1$ or\\
\indent \texttt{sage: E.nloop\_invariant()}\\
otherwise. To do this using a precomputed Neumann-Zagier datum \texttt{Sage}
object file instead of defining \texttt{D} as above use\\
\indent \texttt{sage: nz = load('nz\_exact\_6\_1.sobj')}\\
\indent \texttt{sage: E = nloop(nz, n, diagrams)}.

The entire process described above has been streamlined into two automated
functions for convenience. For example, to start with a specified manifold 
\texttt{M} and diagrams list, compute the Neumann-Zagier datum, and then 
compute the $n$-loop invariant, simply use\\
\indent \texttt{sage: nloop\_from\_manifold(M, n, diagrams, 
engine="retrieve")}.\\
The \texttt{NeumannZagierDatum} optional arguments described above may be
entered here as seen in the example. On the other hand, to start with a
precomputed Neumann-Zagier datum \texttt{Sage} object file (loaded as 
\texttt{nz}) and a diagrams list then compute the $n$-loop invariant, simply 
use\\
\indent 
\texttt{sage: nloop\_from\_nzdatum(nz, n, diagrams,  engine="retrieve")}.

Also available at \cite{nloop-compute} is an almost identical version of our 
code, \texttt{nloop\_num.py}, which produces numerical results to arbitrary 
precision instead of exact computations. The usage for this file is the same.

\subsection{Sample computations}
\lbl{sub.sample}

The results of our computations are available from \cite{nloop-compute},
along with the code and data files.

To illustrate our method, consider the $6_1=\mathrm{K4}_1$ knot with trace
field $F_{6_1}=\BQ(x)$, where $x=-1.50410836415074\dots 
+ i 1.22685163774658\dots $ 
is a root of
$$
x^4 + 2 x^3 + x^2 - 3 x + 1 = 0.
$$
$F_{6_1}$ is a number field of type $[0,2]$ with discriminant $257$, a prime
number. It follows that the Bloch group $\calB(F_{6_1})$ is an abelian group
of rank $2$~\cite{suslin,zi}. The default \texttt{SnapPy} 
triangulation for $\mathrm{K4}_1$ uses 4 ideal tetrahedra with shapes 
\small{
$$
z=\left(\frac{3}{2} x^{3} + \frac{7}{2} x^{2} + 3 x - \frac{5}{2}, \,\,
2 x^{3} + 5 x^{2} + 5 x - 3, \,\,
-\frac{1}{2} x^{3} - \frac{3}{2} x^{2} - x + \frac{3}{2}, \,\,
\frac{1}{2} x^{3} + \frac{3}{2} x^{2} + 2 x + \frac{1}{2} \right).
$$
}
A Neumann-Zagier datum $\ga=(A,B,\nu,f,f'',z)$ is given by
$$
A = \left(\begin{array}{rrrr}
1 & 0 & -1 & 0 \\
-1 & 1 & 1 & 1 \\
-1 & 1 & 1 & 0 \\
-1 & 1 & 0 & 1
\end{array}\right), \,\,
B = \left(\begin{array}{rrrr}
1 & 0 & 0 & 1 \\
0 & 1 & 0 & 1 \\
0 & 0 & 1 & 0 \\
0 & 0 & 0 & 2
\end{array}\right), \,\, 
\nu = \left(\begin{array}{r} 
1 \\ 2 \\ 1 \\ 2 
\end{array}\right), \,\, 
f = \left(\begin{array}{r} 
1 \\ 2 \\ 0 \\ 1
\end{array}\right), \,\, 
f'' = \left(\begin{array}{r} 
0 \\ 0 \\ 0 \\ 0
\end{array}\right).
$$
The $n$-loop invariants for $n=1,\dots,6$ are given by:
{\small
\begin{align*}
\tau &=-\frac{7}{2} x^{3} - \frac{17}{2} x^{2} - \frac{17}{2} x + 6
\\
S_2 &= \frac{46490}{198147} x^{3} + \frac{231209}{396294} x^{2} 
+ \frac{473191}{792588} x - \frac{62777}{264196}
\\
S_3 &= \frac{570416}{16974593} x^{3} + \frac{2833463}{33949186} x^{2} 
+ \frac{1122215}{16974593} x - \frac{1386486}{16974593}
\\
S_4 &= -\frac{2255130587026}{50451970187565} x^{3} 
- \frac{91695358340911}{807231523001040} x^{2} 
- \frac{85651263871967}{807231523001040} x 
+ \frac{1596902056811}{20180788075026}
\\
S_5 &= -\frac{37040877003091}{1728820845093894} x^{3} 
- \frac{330280282463219}{6915283380375576} x^{2} 
- \frac{53499149965837}{1728820845093894} x 
+ \frac{72838757049049}{1152547230062596}
\\
S_6 &= \frac{1449319256564305241317}{17984434859623040256945} x^3 +
\frac{23592842410230239076799}{115100383101587457644448} x^2 +
\\ & \quad\,\,
\frac{110567432832899754708187}{575501915507937288222240} x -
\frac{20008494585620168748319}{143875478876984322055560}
\end{align*}
}

\subsection*{Acknowledgment}
The authors wish to thank the computer support group of the School of 
Mathematics and especially Justin Filoseta and Lew Lefton for their support.
We also wish to thank the anonymous referee for a careful reading of our
article and for suggested improvments in the presentation.

%%%%%%%%%%%%%%%%%%%%%%%%%%%%%%%%%%%%%%%%%%%%%%%%%%%%%%%%%%%%%%%%%%%%%%%%%%%%
%%%%%%%%%%%%%%%%%%%%%%%%%%%%%%%%%%%%%%%%%%%%%%%%%%%%%%%%%%%%%%%%%%%%%%%%%%%%

\bibliographystyle{hamsalpha}
\bibliography{biblio}
\end{document}